\documentclass[12pt]{amsart2000}

\theoremstyle{plain}
\newtheorem{theorem}{Theorem}
\newtheorem{lemma}{Lemma}

\newtheorem{proposition}[theorem]{Proposition}
\theoremstyle{definition}

\newtheorem{remark}{Remark}

\def\Const{{\rm Const\ }}

\def\Lie{\mathop{\rm Lie}}

\def\EXP{{\mathbb{E}}}
\def\PROB{{\mathbb{P}}}

\def\T{{\mathbb{T}}}

\def\gm{\gamma}
\def\eps{\varepsilon}
\def\dt{\delta}
\def\sg{\sigma}

\def\Card{{\#}}

\def\bdef{\begin{Def}}
\def\endef{\end{Def}}
\def\bthm{\begin{theorem}}
\def\ethm{\end{theorem}}
\def\blm{\begin{lemma}}
\def\elm{\end{lemma}}
\def\brm{\begin{remark}}
\def\erm{\end{remark}}
\def\beq{\begin{eqnarray}}
\def\eneq{\end{eqnarray}}
\def\beal{\begin{aligned}}
\def\enal{\end{aligned}}

\def\l{\left}
\def\r{\right}

\def\bL{{\mathbf{L}}}

\def\brC{{\bar C}}
\def\brK{{\bar K}}

\def\brj{{\bar j}}

\def\brgamma{{\bar\gamma}}

\def\brnu{{\bar\nu}}

\def\cB{{\mathcal{B}}}
\def\cC{{\mathcal{C}}}
\def\cF{{\mathcal{F}}}
\def\cG{{\mathcal{G}}}

\def\StB{\mathbb B}
\def\R{\mathbb R}
\def\Z{\mathbb Z}
\def\integers{\mathbb Z}

\def\hj{{\hat j}}

\def\hgamma{{\hat\gamma}}

\def\ta{{\tilde a}}
\def\tb{{\tilde b}}
\def\tPhi{{\tilde \Phi}}

\def\th{\theta}
\def\~{\tilde}
\def\lb{\lambda}

\def\ve{{\vec e}}

\def\RmII{{I\!\!I}}
\def\HD{{\rm HD}}

\author{D. Dolgopyat, V. Kaloshin and L. Koralov}
\title{Hausdorff dimension in stochastic dispersion.}
\address{\begin{tabular}{lll}
Dmitry Dolgopyat & Vadim Kaloshin & Leonid Koralov \\
Department of Mathematics & Department of &
Department of Mathematics \\
PennState University & Mathematics MIT & Princeton University \\
University Park PA 16802 & Cambridge MA 02139 & Princeton NJ 08544 \\
dolgop@math.psu.edu & kaloshin@math.mit.edu &
koralov@math.princeton.edu \\
www.math.psu.edu\slash dolgop\slash  & & \\
\end{tabular}}
\thanks{D.D. was partially supported by NSF and Sloan Foundation,
V. K. was partially supported by American Institute of Mathematics
Fellowship and Courant Institute, and L. K. was partially supported
by NSF postdoctoral fellowship.}
\begin{document}
\begin{abstract}
We consider the evolution of a connected set in Euclidean space carried
by a periodic incompressible stochastic flow. While for almost every
realization of the random flow at time $t$ most of the particles are at
a distance of order $\sqrt{t}$ away from the origin \cite{DKK1},
there is an uncountable set of measure zero of points, which escape
to infinity at the linear rate \cite{CSS1}. In this paper we prove
that this set of linear escape points has full Hausdorff dimension.
\end{abstract}
\maketitle

\centerline{\it Dedicated to our teacher Yakov Sinai
on occasion of his 65th birthday.} \vskip 10mm

\section{Introduction.}

One of the greatest achievements in mathematics of the second half of
the last century was creation of the theory of hyperbolic dynamical
systems in works of Anosov, Bowen, Ruelle, Sinai, Smale and many others.
The importance of this theory is not so much in that it allows one to get
new information about a large class of ordinary differential equations
but in that it provides a paradigm for understanding irregular behavior
in a large class of natural phenomena. From the mathematical point of
view it means that the theory should be useful in many branches of
mathematics beyond the study of finite-dimensional dynamical systems.
The aim of this note is to illustrate this on a simple example. Namely,
we show how the theory of nonuniformly hyperbolic systems, i.e. systems
with non-zero Lyapunov exponents, can explain ballistic behavior
in a problem of passive transport in random media.

This paper concerns the long time behavior of a passive substance
(say an oil spill) carried by a stochastic flow. Various aspects of
such behavior  have been a subject of a number of recent papers (see
\cite{Cm, CC, CGXM, CS, CSS1, CSS2, DKK1, DKK2, LS, SS, ZC1, ZC2} etc.)
Consider an oil spill at the initial time concentrated in a domain
$\Omega$. Let $\Omega$ evolve in time along trajectories of
the stochastic flow and $\Omega_t$ be its image at time $t$.
The papers mentioned study the rate of stretching of the boundary
$\partial \Omega_t$, growth of the diameter and the ``shape'' of
$\Omega_t$, distribution of mass of $\Omega_t$, and many other related
questions. In this paper we model the stochastic flow by a stochastic
differential equation driven by a finite-dimensional Brownian motion
$\{\theta(t)=(\theta_1(t),\dots, \theta_d(t))\in \R^d\}_{t\geq 0}$
\begin{equation}
\label{SF}
dx_t=X_0(x_t)dt+\sum_{k=1}^d X_k(x_t) \circ d\theta_k(t)
\quad x\in \R^N
\end{equation}
where $\{X_k\}_{k=0}^d$ are $C^\infty$-smooth space periodic divergence
free vector fields on $\R^N$. Alternatively one can regard this system
as a flow  on $\T^N=\R^N/\integers^N$.\
Below we impose certain nondegeneracy assumptions on vector fields
$\{X_k\}_{k=1}^d$ from \cite{DKK1}. These assumptions hold
on an dense open set of $C^\infty$-smooth divergence free vector
fields on $\T^N$ or satisfied generically.

An interesting feature of the flow (\ref{SF}) is the dichotomy
between growth of the mass and shape of the spill $\Omega_t$. On one
hand, most of the points of the tracer $\Omega_t$ move at distance
of order $\sqrt{t}$ at time $t$. More precisely, let $\rho$ be
a smooth metric on $\T^N$, naturally lifted to $\R^N$ and $\nu$
be a measure of a finite energy, i.e. for some positive $p$ we have
$$
\int\!\!\!\!\int\frac{d\nu(x)d\nu(y)}{\rho^p(x,y)}<\infty.
$$
In particular, $\nu$ can be the Lebesgue probability measure supported
on an open set $\Omega$, which also supports the initial oil spill.
Let $\nu_t$ be its image under the flow (\ref{SF}) and $\brnu_t$ be
rescaling of $\nu_t$, defined by as follows: for a Borel set
$\Omega\subset\R^N$ put $\brnu_t(\Omega)=\nu_t(\sqrt{t} \Omega)$.

\begin{theorem} \label{ThMass} (\cite{DKK1}) For almost every
realization of the Brownian motion $\{\theta(t)\}_{t\geq 0}$
the measure $\brnu_t$ weakly converges to a Gaussian measure
on $\R^N$ as $t\to\infty$.
\end{theorem}

\brm
Notice that this is the Central Limit Theorem with respect to
randomness in initial conditions, not with respect to randomness
of the Brownian motion $\{\theta(t)\}_{t\geq 0}$.
\erm

On the other hand, there are many points with {\it linear growth.}
Fix a realization of the Brownian motion $\{\th(t)\}_{t\geq 0}$.
Let $\bL_\theta$ denote the random set of points with a linear escape
$$
\bL_\th=\left\{x\in \R^N:
\ \lim\inf_{t\to+\infty}\frac{|x_t|}{t}>0\right\}.
$$
The following result is a special case of \cite{SS}
(see also \cite{CSS1}).

\begin{theorem}
\label{ThLin}
Let $S$ be a connected set containing at least two points.
Then for almost every realization of the Brownian motion
$\{\theta(t)\}_{t\geq 0}$ the set $\bL_\theta \bigcap S$ is
uncountable.
\end{theorem}

In fact in dimension 2 there is a limiting shape of the rescaled
contaminated area. Namely, let $\Omega$ be a bounded open set,
$\Omega_t$ be its image under the flow (\ref{SF}) and
$\cC_t= \bigcup_{0\leq s\leq t} \Omega_s.$ In other words,
we call a point $x$ {\it contaminated by the time $t$} if there
is a trajectory from our curve which has passed through $x$
before time $t$.

{\bf The Shape Theorem.} {\it (\cite{DKK2})
If $N=2$, then there exists a convex compact set $\StB\subset \R^2$
such that for almost every realization of the Brownian motion
$\{\theta(t)\}_{t\geq 0}$ and any $\delta>0$ there
exists $T=T(\delta)$ such that for all $t>T$
$$
(1-\delta)\ \StB \subset \frac{\cC_t}{t}
\subset (1+\delta)\ \StB.
$$
}
\begin{remark}
The ``shape'' $\StB\subset \R^2$ is independent of the initial
spill $\Omega$. Moreover, an open set $\Omega$ can be replaced by
a smooth curve $\gamma$ for the Shape Theorem to hold true.
\end{remark}
In view of Theorems \ref{ThMass}, \ref{ThLin}, and  the Shape
Theorem it is interesting  to see how large is the set of points
with linear growth. In this paper we first prove the following

\begin{theorem} \label{HDC}
Let $\gamma$ be a smooth curve on $\R^2$. Then for almost every
realization of the Brownian motion $\{\theta(t)\}_{t\geq 0}$ we have
$\HD(\bL_\th\bigcap\gamma)=1.$
\end{theorem}

Then in Section \ref{fullHD} using this Theorem we derive
the following main result of the paper
\begin{theorem} \label{FullDim} ({\bf Main Result}) For almost every
realization of the Brownian motion $\{\theta(t)\}_{t\geq 0}$ we have
that points of the flow (\ref{SF}) with linear escape to infinity $\bL_\th$
form a dense set of full Hausdorff dimension $\HD(\bL_\theta)=N.$
\end{theorem}

By Theorem \ref{ThMass} for most points $x_0=x$ in $\R^N$
its trajectory $x_t$ is of order $~\sqrt{t}$ away from the origin
at time $t$. Also, the Law of Iterated Logarithm for functionals of
diffusion processes and Fubini Theorem imply that
the set of points $\bL_\th$ with linear
escape has measure zero. This Corollary says that $\bL_\th$ is
the ``richest'' possible set of measure zero in $\R^N$, namely, is
of full Hausdorff dimension $N$.

\section{Nondegeneracy assumptions.}
\label{asm}  In this section we formulate a set of assumptions on
the vector fields, which in particular imply the Central Limit Theorem
for measures, the estimates on the behavior of the characteristic
function of a measure carried by the flow (see \cite{DKK1}), and large
deviations estimates (see \cite{BS}). Such estimates are essential for
the proof of our results. Recall that $X_0, X_1, \dots , X_d$ are
assumed to be $C^\infty$-smooth, periodic and divergence free.

(A) ({\it hypoellipticity for $x_t$}) For all $x\in \R^N$ we have
\[
{\Lie} (X_1,\dots,X_d)(x)=\R^N.
\]

Denote the diagonal in $\T^N \times \T^N$ by
\[
\Delta=\{(x^1,x^2)\in \R^N \times \R^N:\ x^1=x^2\
({\rm mod}~ 1) \}.
\]

(B) ({\it hypoellipticity for the two--point motion}) The
generator of the two--point motion $\{(x^1_t,x^2_t):\ t>0\}$ is
nondegenerate away from the diagonal $\Delta$, meaning that the
Lie brackets made out of $(X_1(x^1),X_1(x^2)), \dots,
(X_d(x^1),X_d(x^2))$ generate $\R^N \times \R^N$.

To formulate the next assumption we need additional notations. Let
$Dx_t:T_{x_0}\R^N \to T_{x_t}\R^N$ be the linearization of $x_t$
at $t$. We need the hypoellipticity  of the process
$\{(x_t,Dx_t):\ t>0\}$. Denote by $TX_k$ the derivative of the
vector field $X_k$ thought as the map on $T\R^2$ and by
$S\R^N=\{v\in T\R^N:\ |v|=1\}$ the unit tangent bundle on $\R^N$.
If we denote by $\~X_k(v)$ the projection of $TX_k(v)$ onto
$T_vS\R^N$, then the stochastic flow (\ref{SF}) on $\R^N$ induces
a stochastic flow on the unit tangent bundle $S\R^N$, defined by
the following equation:
 \[ d\~x_t=\sum_{k=1}^d \~ X_k(\~x_t)\circ
d\th_k(t)+ \~X_0(\~x_t)dt. \]
With these notations we have
condition

(C) ({\it hypoellipticity for $(x_t,Dx_t)$}) For all $v\in S\R^N$
we have
\[
\Lie(\~X_1,\dots,\~X_d)(v)=T_vS\R^N~.
\]

For measure-preserving stochastic flows with conditions (C)
Lyapunov exponents $\lb_1, \dots , \lb_N$  exist by {\it multiplicative
ergodic theorem for stochastic flows} of diffeomorphisms (see
\cite{Cv}, Thm. 2.1). Moreover, the sum of Lyapunov exponents
$\sum_{j=1}^N \lb_j$ should be zero (see e.g. \cite{BS}). Under
conditions (A)--(C) the leading Lyapunov exponent is positive
\beq \label{positive}
\lb_1=\lim_{t\to \infty}
\frac {\log |d\varphi_t(x)(v)|}{t} > 0~,
\eneq
where $d\varphi_t(x)$ is the linearization matrix of
the flow (\ref{SF}) integrated from $0$ to $t$ at the point $x$.
Indeed, Theorem 6.8 of \cite{Bx} states that under condition (A)
the maximal Lyapunov exponent $\lb_1$ can be zero only if for
almost every realization of the flow (\ref{SF}) one of the following
two conditions is satisfied

(a) there is a Riemannian metric $\rho'$ on $\T^N$, invariant with
respect to the flow (\ref{SF}) or

(b) there is a direction field $v(x)$ on $\T^N$ invariant with
respect to the flow (\ref{SF}).

However (a) contradicts condition (B). Indeed, (a) implies that
all the Lie brackets of $\{(X_k(x^1),\ X_k(x^2))\}_{k=1}^d$
are tangent to the leaves of the foliation
\[
\{(x^1,x^2)\in \T^N \times \T^N:\ \rho'(x^1,x^2)=\Const\}
\]
and don't form the whole tangent space. On the other hand (b)
contradicts condition (C), since (b) implies that all the Lie
brackets  are tangent to the graph of $v$. This positivity of
$\lb_1$ is crucial for our approach.
\begin{remark}Let us mention an important difference between deterministic
and stochastic dynamics. Most of the results dealing with statistical
properties of deterministic systems assume that all Lyapunov exponents
are non-zero. By contrast we need only one positive exponent. This is
because in the random situation hypoellipticity condition (C) implies that
growth rate of any deterministic vector is given by the largest exponent
(see equation (\ref{GrVec}). This allows us to get our results without
assuming that all the exponents are non-zero.
\end{remark}

We further require that the flow has no deterministic drift, which
is expressed by the following condition

(E) ({\it zero drift})
\[
\int_{\T^2} \l(\sum_{k=1}^d L_{X_k} X_k + X_0\r)(x) dx =0~,
 \]
where $L_{X_k} X_k(x)$ is the derivative of $X_k$ along $X_k$ at
the point $x$. Notice that $\sum_{k=1}^d L_{X_k} X_k + X_0$
is the deterministic components of the stochastic flow (\ref{SF})
rewritten in Ito's form.

The Central Limit Theorem for measures was formulated in \cite{DKK1}
under an additional assumption
\beq \label{zd}
\int_{\T^2} X_k(x) dx =0~,~~k=1,...,d~.
\eneq
This assumption is not needed for the proof of Theorem \ref{HDC}
and as the result for the proof of the Main Theorem. However, in order
to simplify the proof, i.e. use the results of \cite{DKK1} without
technical modifications, we shall assume (\ref{zd}) to hold.

\section{Idea of the proof.}

\subsection{A Model Example}\label{model}

Below we define a random dynamical system on $\R$ which models
the motion of the projection of the spill $\Omega_t$ onto
a fixed line $l \subset \R^N$.

Introduce notations: $I(b;a)=[b-a/2,b+a/2]$ --- the segment
on $\R$ centered at $b$ of length $a$; $s \in \{0,1\}^{\Z_+}$
a semiinfinite sequence of $0$'s and $1$'s, $s_k \in \{0,1\}^k$
a set of $k$ numbers $0$ or $1$,
$\{\{\theta_{s_k}(t)\}_{s_k\in \{0,1\}^k}\}_{k\in\Z_+}$
countable number of standard i.i.d. Brownian motions on $\R$
indexed by binary sequences. Let $\tau$ be positive.

The random dynamical system is defined as follows. Let
$I^\emptyset=I(0;1)$. Then $\sg^\th_0: I^\emptyset \to \R$
stretches $I^\emptyset$ uniformly by $2$ around its center and
shifts it randomly by $\th_\emptyset(\tau)$. Divide
$\sg^\th_0(I^\emptyset)$ in two equal parts $I^0$ and $I^1$
\beq \label{rds1}
\sg^\th_0(I^\emptyset)=I^0 \cup I^1 =
I(\th_\emptyset(\tau)-1/2;1) \cup I(\th_\emptyset(\tau)+1/2;1).
\eneq
Now $\sg^{\th}_1$ acts on each $\{I^i\}_{i=0,1}$ independently by
stretching each $I^i$'s uniformly by $2$ around its center and
shifting by $\th_0(\tau)$ and $\th_1(\tau)$ respectively.
\beq \label{rds2}
\beal
\sg^\th_1 \circ \sg^\th_0 (I^\emptyset)= (I^{00} \cap I^{01}) &
\cup (I^{10} \cup I^{11}) = \\
 I([\th_\emptyset(\tau)-1/2]+[\th_0(\tau)-1/2];1) & \cup
 I([\th_\emptyset(\tau)-1/2]+[\th_1(\tau)+1/2];1)
\\ \cup\ \
I([\th_\emptyset(\tau)+1/2]+[\th_0(\tau)-1/2];1) & \cup
I([\th_\emptyset(\tau)+1/2]+[\th_1(\tau)+1/2];1),
\enal
\eneq
and so on.

Let $n \in \Z_+$. Then at the $n$-th stage ``after time $n\tau$''
the image of the initial unit interval $I^\emptyset=[-1/2,1/2]$
consists of $2^n$ unit intervals. The preimage of each of those unit
intervals is an interval of length $2^{-n}$ uniformly contracted.
Let's give a different definition of the random dynamical system
(\ref{rds1})-(\ref{rds2}).

Consider an isomorphism of the dynamical system on the unit
interval $I=I^\emptyset+1/2=[0,1]$ given by
$\phi:x\mapsto 2x\ (\textup{mod}\ 1)$ and the one sided Bernoulli
shift on two symbols, say $0$ and $1$. Such an isomorphism is given
by $s:x \mapsto s(x)=\{s_k(x)\}_{k=0}^\infty
\in \{0,1\}^{\Z_+}$, where for each $k\in \Z_+$
\beq
\begin{cases}
s_k(x)=0 \quad \textup{if} \ \ \ \phi^n(x)<1/2 \\
s_k(x)=1 \quad \textup{otherwise}.
\end{cases}
\eneq
Let $\eta_n(x)=\#\{k\leq n: s_k(x)=1\}$.
Notice now that
\beq
\sg^\th_n \circ \sg^\th_{n-1} \circ \cdots \circ \sg^\th_0 (x)=
\sum_{k=0}^n \th_{s_k}(\tau)+(\eta_{n+1}(x)-(n+1)/2)/2,
\eneq
where $\th_{s_k}(\tau)$'s are i.i.d. Brownian motions.
Define $\eta_-(x)=\liminf_{n\to\infty} \eta_n(x)/n$. Then for
almost all points $x \in I$ we have
$\eta_-(x)=\lim_{n\to\infty} \eta_n(x)/n = 1/2$.
Let us show however that there is full Hausdorff dimension
set of points in  the interval $I$ such that frequency of $0$'s
is less than frequency of $1$'s, i.e. $\HD\{x\in I:\eta(x)>1/2\}=1$.
Since $\sum_{k=0}^n \th_{s_k}(\tau)/n \to 0$ almost surely
this would imply that the set of points in $I^\emptyset$ with
a nonzero drift for the random dynamical system, defined by
(\ref{rds1})-(\ref{rds2}), has full Hausdorff dimension almost
surely, but is of measure zero.

We shall justify the fact that $\HD\{x\in I:\eta(x)>1/2\}=1$.

\subsection{Points with a nonzero drift} \label{drift}

Fix an arbitrary small positive $\eps$. The goal is to find a fractal
set of points $I_\infty\subset I^\emptyset$ and a probability measure
$\mu_\infty$ supported on $I_\infty$ such that $\mu_\infty$-a.e.
point $x\in I_\infty$ has a nonzero drift to the right, i.e.
$\liminf_{n\to \infty} \sg^\th_n \circ \dots \circ \sg^\th_0 (x)/n>0$.
Moreover, $\HD(\mu_\infty)$ tends to $1$ as $\eps$ tends to $0$.

Construction of the set $I_\infty$ and of the measure $\mu_\infty$ is
inductive. $I_\infty$ is defined as a countable intersection of
a nested sequence of compact sets and $\mu_\infty$ is given as a weak
limit of Lebesgue measures supported on those sets. We describe
the base of the induction and the inductive steps.

$\bullet$ For $n=1$ we have $\sg^\th_0(I^\emptyset)$ is a segment of
length $2$ or union of two segments $I^0$ and $I^1$ of length $1$
each. Cut off the bottom $\eps$-segment from each segment.
This corresponds to cutting off $\eps$-segments $[-1/2,-1/2+\eps]$
and $[0,\eps]$ from $I^\emptyset$.  Denote the surgery result by
$I_0\subset I^\emptyset$ and by $\mu_0$ the Lebesgue probability
measure supported on whole $I_0$. Notice that
\beq \label{drifter}
\mu_0\{x\in I_0: \sg^\th_1(x)=1\}>
\mu_0\{x\in I_0: \sg^\th_1(x)=0\}
\eneq
creates a nonzero drift up, since frequency of $1$'s exceeds
frequency of $0$'s.

$\bullet$ Suppose $I_{n-1}$ and $\mu_{n-1}$ are constructed.
To construct $I_n$ and $\mu_n$ consider the image $\sg^\th_n(I_{n-1})$.
It consists of $2^n$ segments of equal length close to $1$. Cut off
the bottom $\eps$-segment from each. This corresponds to cutting off
$2^n$ segments of length $2^{-n}\eps$ from $I_{n-1}\subset I^\emptyset$.
The result of the surgery is denoted by $I_n$ and by $\mu_n$ we denote
the Lebesgue probability measure supported on the whole $I_n$. Again
the surgery increases probability of $s_n(x)$ being $1$ over $s_n(x)$
being $0$. Thus, this creates a positive drift.

The intersection $I_\infty = \cap_n I_n$ is a fractal set
and the weak limit measure $\mu_\infty=\lim \mu_n$ has Hausdorff
dimension approaching $1$ as $\eps$ tends to zero.
It follows from the construction that for $\mu_\infty$-almost every
point $\eta_-(x)=\mu_0\{x\in I_0: \sg^\th_1(x)=1\}>1/2$.

\subsection{Difficulties in extending of the Model Example
to the case of the flow}

Let $\gm \subset \R^N$ be a smooth curve, $l \in \R^N$ be a line,
and $\pi_l:\R^N \to l$ be an orthogonal projection onto $l$. Suppose
at the initial moment of time $\pi_l(\gm)=I^\emptyset$ is the unit
interval. If not, then rescale it and shift it to the origin.

The most subtle element in extending  the Model Example is defining
the stopping (stretching) time $\tau$ or deciding {\it when to stop}
$\gm_t$ and {\it how to cut off} some parts of $\gm_t$ in order to create
a nonzero drift as in (\ref{drifter}). Such a stopping time needs
to have several important features\footnote{In the Model Example
$\tau$ is a constant}.

1. {\it Stretching property of the stopping time}: It is not difficult
to show that if $|\pi_l(\gm_0)|=1$, then the stopping time
\beq
\tau_\gm=\inf \{ t \geq 0: |\pi_l(\gm_t)|=2\}
\eneq
has finite expectation and exponential moments uniformly bounded over
all compact curves with projection of length 1 (see e.g. \cite{CSS1}).

The analogy between this $\tau_\gm$ and the model $\tau$ is clear.
However, the geometry of $\gm_{\tau_\gm}$ in $\R^N$ might become quite
complicated ($\gm_{\tau_\gm}$ might spin, bend, fold, and so on see
computer simulations in \cite{CC}) so it is not reasonable to stop
all parts of the curve $\gm$ simultaneously and perform the surgery
(cut off of ``bottom'' parts as in the passage preceding (\ref{drifter})).
For this reason at the first stage of a partition/cut off process
we split $\gm_{\tau_\gm}$
not in {\it two} parts as in the Model Example, but in {\it a number}
(may be countable) of random parts $\gm=\cup_{j\in J}\gm_j$ and each
part $\gm_j$ will have {\it its own stopping time} $\tau_j$.

2. {\it A countable Partition of $\gm$}: We shall partition $\gm$ into
at most countable number of segments $\gm=\cup_{j\in J}\gm_j$ (see
Section \ref{ScPart}). Each $\gm_j$ has its own stopping time $\tau_j$
so that the image $\varphi_{\tau_j} \gm_j$ under the flow (\ref{SF})
is not too folded (see condition (b) of Theorem \ref{ThHypT}).
Moreover, such a stopping time $\tau_j$ still has finite expectation
and exponential moments (see condition (e) of Theorem \ref{ThHypT}).

Now if we have that the image $\varphi_{\tau_j}\gm_j$ is ``regular''
it {\it does not reflect dynamics} on $\gm_j$. In order
to imitate the Model Example's cut off construction we need to stop
$\varphi_t\gm_j$ at  the moment when $\varphi_t|_{\gm_j}$ is more or
less uniformly expanding on $\gm_j$ forward in time or
$(\varphi_t)^{-1}|_{\varphi_t\gm}$ is uniformly contracting on
$\varphi_t \gm$ backward in time (see conditions (a), (c),
and (d) of Theorem \ref{ThHypT}). For example, if there is no backward
contraction by dynamics of $\varphi^{-1}_{\tau_j}$ on
$\varphi_{\tau_j}\gm_j$, then if we cut off an $\eps$-part of
$\varphi_{\tau_j}\gm_j$ its preimage in $\gm_j$ might not be small
compare to length of $\gm_j$. As we explain in more details below
we need this smallness to estimate Hausdorff dimension of remaining
points in $\gm\supset \gm_j$ after the surgery. The property of
uniformness of distortion of a dynamical system is usually called:

3. {\it A bounded distortion property}: In the Model Example, Section
\ref{model} we have {\it uniform} backward contraction of intervals:
at stage $n$ (after time $\tau n$) by a factor $2^{-n}$. So, when
we cut off an $\eps$-part of an interval at stage $n$, it
corresponds to $2^{-n}\eps$-part of the initial segment $I^\emptyset$.
This remark makes an estimate of Hausdorff dimension of the set
$I_\infty$ or of the measure $\mu_\infty$ supported on $I_\infty$
trivial, because the sets $\{I_n\}_{n\in \Z_+}$  have selfsimilar
structure. Certainly, this is no longer true for the evolution of
$\gm$ under the flow (\ref{SF}). Some parts of $\gm$ expanded by
$\varphi_t|_{\gm}$ expanded more than others. Condition (c) of
Theorem \ref{ThHypT} makes sure that there is a backward contraction
in time and condition (d) of the same theorem says that rate of
backward contraction Holder regularly depends on a point on a short
interval $\gm_j$. Thus, backward contraction is sufficiently uniform
on $\gm_j$'s.

In the theory of deterministic dynamical systems with non-zero Lyapunov
exponents the set of points satisfying uniform estimates for forward
and backward expansion (as well as uniform estimates for angles between
stable and unstable manifolds) are called {\it Pesin sets} and times
when an orbit visits given Pesin set are called {\it hyperbolic times}.
The existence of Pesin sets follows from abstract
ergodic theory (see \cite{P1}). Understanding the geometry of
these sets in concrete examples is an important but often difficult
task. In this paper we describe some properties of Pesin sets for
stochastic flows. This description plays a key role in the proof of
Theorem \ref{HDC} and we also think it can be useful in many other
questions about stochastic flows.

In particular let us mention that the estimates similar to ones given
in Section \ref{ScHypM} play important role in many other questions in
the theory of deterministic systems such as periodic orbit estimates
\cite{K} and constructions of maximal measures \cite{N}, etc.

Our arguments in this paper are quite similar to \cite{D1}, \cite{KM}
even though the control of the geometry of images of curves is much
more complicated in our case. Some interesting formulas for dimensions
of nontypical points can be found in \cite{BaSc}. We also refer
the readers to the survey \cite{Sz} and the book \cite{P3} for more
results about dimensions of dynamically defined sets.

The rest of the paper is organized as follows. In Section \ref{ScHypM}
in Theorem \ref{ThHypT} we define a stopping time $\tau$ and prove that
it has finite expectation and exponential moments. In
Section \ref{hypmoments} we investigate expansion properties of
the flow (\ref{SF}) at the stopping time $\tau$ and complete the proof
of Theorem \ref{ThHypT}. Recall that section \ref{drift} above was
devoted to the construction of points with a nonzero drift. Namely,
we need to construct a Cantor set $I_\infty$ and
a measure $\mu_\infty$ supported on $I_\infty$ so that $\mu_\infty$-a.e.
point has a nonzero drift. First, in Section \ref{ScPart} we present
an algorithm of construction of a random Cantor set $I$ inside
the initial curve $\gm$. Then, in Section \ref{ScDefMes} we define
a probability measure $\mu$ supported on $I$ with almost sure nonzero
drift. Hausdorff dimension of such a measure is estimated in
Section \ref{ScHD}. Main Result (Theorem \ref{FullDim}) is derived
from Theorem \ref{HDC} in Section \ref{fullHD}. Auxiliary lemmas are
in the Appendix at the end of the paper.

\section{Hyperbolic moments. Control of the smoothness.}
\label{ScHypM}

Introduce notations. Denote by $\varphi_{t_1, t_2}$ a diffeomorphism
of $\T^N$, obtained by solving (\ref{SF}) on the time interval
$[t_1, t_2]$, and by $\varphi_t$ the diffeomorphism  $\varphi_{0,t}.$
The flow (\ref{SF}) can also be thought as the product of independent
diffeomorphisms $\{\varphi_{n, n+1}:\T^N \to \T^N\}_{n\in\Z_+}$.

Given positive numbers $K$ and $\alpha$ we say that a curve $\gamma$
is {\it $(K,\alpha)$-smooth} if in the arclength parameterization
the following inequality holds
\beq  \nonumber
\left| \frac{d\gamma}{ds}(s_1)-\frac{d\gamma}{ds}(s_2)\right|\leq
K \rho^\alpha(s_2,s_1) \ \  \textup{for each pair of points}
\ \  s_1,\ s_2 \in \gamma.
\eneq
In all the inequalities which appear below the distance $\rho$ between
the points on $\gamma$ or its images $\varphi_t\gamma$'s is measured
in the arclength metric induced on $\gamma$ or $\varphi_t\gamma$ from
the ambient space. In order to do not overload notations we omit
dependence on $\gamma$ or $\varphi_t\gamma$ when it is clear from
the context which curve we use.

The goal of this section is to show that for a sufficiently small
$\alpha$ and a sufficiently large $K,$ starting with an arbitrary
point $x$ on a $(K, \alpha)$-smooth curve $\gm$, the part of image
of this curve in a small neighborhood of the image of $x$ is
often smooth. More precisely, we prove the following statement.
Let $\lb_1$ be the largest Lyapunov exponent of the flow (\ref{SF})
which is positive see (\ref{positive}).

\begin{theorem} \label{ThHypT} For any $0<\lb_1'<\lb_1$ there
exist sufficiently small $r>0,\;\alpha\in(0,1)$,
and sufficiently large $K>0$ and $n_0\in \Z_+$ with
the following properties:

For any $(K,\alpha)$-smooth $\gamma$ of length between
$\frac{r}{100}$ and $100 r$ and each point $x\in\gamma$ there is
a stopping time $\tau=\tau(x)$, divisible by $n_0$, such that

\noindent\ (a) $\|d\varphi_\tau|T\gamma(x)\|>100$ and length
of the corresponding curve $l(\varphi_\tau\gamma)\geq r;$

Denote by $\brgamma_r$ a curve inside $\varphi_\tau\gamma$ of
radius $r$ with respect to induced in $\varphi_\tau\gamma$ length
centered at $\varphi_\tau(x).$ Then

\noindent (b)\ $\brgamma_r$ is $(K,\alpha)$-smooth

\noindent and for each pair of points $y_1, y_2 \in \brgamma_r$
the following holds

\noindent (c)\ for each integer $0\leq k\leq \frac{\tau}{n_0}$
we have
$$
\rho(\varphi_{\tau, \tau-k n_0} y_1,\varphi_{\tau, \tau-k n_0} y_2)
\leq e^{-\lb_1' k n_0} \rho(y_1, y_2);
$$

\noindent (d)\ $\left|  \ln \|d\varphi_\tau^{-1}|T\brgamma_r\|(y_1)-
\ln \|d\varphi_\tau^{-1}|T\brgamma_r\|(y_2) \right|
\leq\Const \rho^\alpha(y_1, y_2)$;

Moreover, for such a stopping time $\tau(x)$ we have

\noindent (e)\ $\EXP\, \tau(x)\leq C_0;\
\PROB\{\tau(x)>T\}\leq C_1 e^{-C_2 T}$ for any $T>0$,

All the above constants depend only on vector fields $\{X_k\}_{k=0}^d$
and $\lb_1$, but independent of the curve $\gm$.
\end{theorem}

\brm Choosing integer $n_0$ is only for our convenience. Requirement
that $\tau$ is divisible by $n_0$ will be used for construction of
partition of $\gm$ in Section  \ref{ScPart}.
This is also indication of flexibility in choice of both constants.
The choice of constants 100, 1000, etc. in this paper is more or
less arbitrary. Any constant greater than 1 would suffice.
\erm

\proof The leading idea of the proof is that with probability close to 1
for a sufficiently large $n_0$ the diffeomorphism
$\varphi_{t,t+n_0}:\T^N\to \T^N$ gets close to its asymptotic behavior.
In particular, the norm of the linearization $\|d\varphi_{t,t+n_0}(x)\|$
as the matrix is $\sim \exp(\lb_1n_0+o(n_0))$ as the top Lyapunov exponent
predicts. Moreover, the linearization dominates higher order terms of
$\varphi_{t,t+n_0}(x)$ and, therefore, determines local dynamics in
a neighborhood of $x$. Thus, to some extend for large periods of time
the flow (\ref{SF}) behaves similarly to uniformly hyperbolic system,
for which properties of the Theorem are easy to verify. Now we start
the proof.

First we construct a stopping time $\tau$ as a first moments
satisfying a certain number of regularity inequalities
(see (\ref{ScDer})--(\ref{FDer})).  This inequalities would include
 $K,\ r,\ n_0$ and some other parameters. Then we show that for any
$\varepsilon>0$ these parameters can be adjusted so that probability
that the number of times each inequality is violated up to time $T$
at least $\varepsilon T$ times decay exponentially in $T.$ This would
guarantee condition (e). Finally, we show that these inequalities
imply conditions (a)--(d) and as the result prove the Theorem. In
the Appendix we obtain large deviation estimates necessary for the
proof below.

Our first goal is to control distortion of the unit tangent vector to
images $\varphi_t \gm$ of $\gamma$ as time $t$ evolves.
Consider a collection of subsets of $\gamma$ indexed by $j$
$$
\cB_{T,j n_0}(x)=\{y\in\gamma:
\rho(\varphi_{n_0 k}y, \varphi_{n_0 k}x)\leq r e^{-\lb_1'(T-n_0 k)}\;
\text{for} \; 0\leq k\leq j\},
$$
where $j$ varies from $1$ to $T/n_0$. We would like to find
an integer moment of time $\tau$, divisible by $n_0$, such that
\begin{equation}
\label{SlGrCu}
\varphi_{j n_0} \cB_{\tau,\tau}(x) \; \text{is}\;
(K e^{\epsilon(\tau-j n_0)},\alpha)- \text{smooth for all}\
\ \ j=0,\dots ,\frac{\tau}{n_0}.
\end{equation}
The rest of this section is devoted to showing that the set of those $T$,
divisible by $n_0$ for which (\ref{SlGrCu}) holds with $T=\tau$ has
density close to 1 if $K$ is sufficiently large. Given $T$ denote by
$K_j$ the $\alpha$-Holder norm of $\varphi_{j n_0} \cB_{T,j n_0}(x).$
We would like to derive an inductive in $j$ formula relating $K_j$ and
$K_{j+1}$ so that in $T/n_0$ steps we get a required statement. Let
$z_1,z_2$ be two points on $\varphi_{j n_0} \cB_{T,j n_0}(x)$ and $r$
is sufficiently small, then
$$
\rho(\varphi_{j n_0, (j+1) n_0} z_1, \varphi_{j n_0, (j+1) n_0} z_2)\geq
\frac 12 \inf_{\varphi_{j n_0}\cB_{T,j n_0}(x)}
\|d\varphi_{j n_0, (j+1) n_0}| T\gamma\|  \ \rho(z_1, z_2).
$$
Let $d^2 \varphi_{j n_0, (j+1) n_0}$ be the Hessian matrix consisting of
second derivatives of the diffeomorphism $\varphi_{j n_0, (j+1) n_0}$.
Assuming now that for each integer $j< \frac{T}{n_0}$ and some $R>0$
we have
\begin{equation}
\label{ScDer}
\|d^2 \varphi_{j n_0, (j+1) n_0}\|\leq R e^{\epsilon(T-j n_0)}
\end{equation}
and that condition (\ref{SlGrCu}) holds true up for each $j\leq j^*$.
Then we get
\beq
\beal
\rho(\varphi_{j n_0,(j+1) n_0}z_1,\varphi_{j n_0,(j+1) n_0} z_2)
\geq & \
\\ \frac 12
\left(\|d\varphi_{j n_0, (j+1) n_0}| T\gamma\|(\varphi_{j n_0} x)-
rR K e^{-(\lb_1' n_0-2\epsilon)(T-j^*)}\right) & \ \rho(z_1,z_2).
\enal
\eneq
We would like to prove that (\ref{SlGrCu}) holds true for $j=j^*+1$.
Assume also that for each $j<T/n_0$ we have
\begin{equation}
\label{SlDecDr}
rR K e^{-(\lb_1'n_0-2\epsilon)(T-j n_0)} \leq
\frac{\|d\varphi_{j n_0, (j+1) n_0}| T\gamma\|(\varphi_{j n_0} x)}{4}
\end{equation}
then we get
$$
\rho(\varphi_{j^* n_0, (j^*+1) n_0} z_1, \varphi_{j^* n_0, (j^*+1) n_0} z_2)
\geq
\frac{\|d\varphi_{j^* n_0, (j^*+1) n_0}| T\gamma\|(\varphi_{j^* n_0} x)}{4}
\ \rho(z_1, z_2).
$$
Let $v_1$ and $v_2$ be directions of the tangent vectors to
$\varphi_{j^* n_0} \gamma$ at $z_1$ and $z_2$ respectively, then
\beq \nonumber
\beal
\rho(d\varphi_{j^* n_0, (j^*+1) n_0} (z_1) v_1,
d\varphi_{j^* n_0, (j^*+1) n_0} (z_2)v_2)
\leq \\
\rho(d\varphi_{j^* n_0, (j^*+1) n_0} (z_1) v_1,
d\varphi_{j^* n_0, (j^*+1) n_0} (z_1)v_2)  \\
+  \rho(d\varphi_{j^* n_0, (j^*+1) n_0} (z_1) v_2,
d\varphi_{j^* n_0, (j^*+1) n_0} (z_2)v_2).
\enal
\eneq
Denote the first and the second terms by $I$ and $\RmII$
respectively. If (\ref{ScDer}) holds, then
$$
\RmII\leq R e^{\epsilon(T-j^*n_0)} \rho(z_1, z_2).
$$
Now since $v_1$ and $v_2$ are close $T(\varphi_{j^*n_0}\gamma)(x)$
and $z_1$ is close to $\varphi_{j^* n_0} x$ we have
$$
\left[d\varphi_{j^* n_0, (j^*+1) n_0} (z_1) v_1-
d\varphi_{j^* n_0, (j^*+1) n_0} (z_1)v_2 \right]
\approx d_v d\varphi_{j^*n_0, (j^*+1) n_0} (x) (v_1-v_2).
$$
Let us give more precise estimates. Notice that if $A$ is a linear
map, then its action on the projective space satisfies
$$
\|dA(v)\delta v\|=\frac{\|\Pi_{(Av)^\bot} A\delta v\|}{\|Av\|}
\leq \frac{\|A\|}{\|Av\|},
$$
where $\Pi_{(Av)^\bot}$ is the orthogonal projection onto
the direction $Av$ and $\delta v$ is an element of $T_vT_xM$.

Therefore, we can assume that for a positive integer $T/n_0$ and each
$j^*<T/n_0$ the following inequality is satisfied
\begin{equation}
\label{CloseToDer}
\frac{\rho(Av_1, A v_2)}{\rho(v_1, v_2)}\leq
2\frac{\|d\varphi_{j^*n_0, (j^*+1) n_0} \|}
{\|d\varphi_{j^*n_0, (j^*+1) n_0}|T\varphi_{j^*n_0} \gamma\|},
\end{equation}
for any linear map $A$ such that
$\| A-d\varphi_{j^*n_0, (j^*+1) n_0}\| \leq
r R e^{(\lb_1' n_0-\epsilon) (T-j^*n_0)}$ and for any pair $(v_1, v_2)$
of tangent vectors such that
$\rho(v_1, v_2)\leq K e^{-(\lb_1' n_0 \alpha-\epsilon)(T-j^*n_0)}.$

Thus
$$
I\leq 2 K_{j^*}\ \rho^\alpha (z_1, z_2)\
\frac{\|d\varphi_{j^*n_0, (j^*+1) n_0} \|}
{\| d\varphi_{j^*n_0, (j^*+1) n_0}|T\varphi_{j^*n_0} \gamma \|}.
$$
Hence (\ref{ScDer})--(\ref{CloseToDer}) imply that
$$
K_{j^*+1}\leq \frac{8 K_{j^*}\ \|d\varphi_{j^*n_0, (j^*+1) n_0}\|}
{\|d\varphi_{j^*n_0, (j^*+1) n_0}|T\varphi_{j^* n_0}\gamma\|^{1+\alpha}}+
\frac{2 R e^{-(\lb_1' n_0(1-\alpha)-\epsilon)(T-j^*n_0)}}
{\|d\varphi_{j^* n_0, (j^*+1) n_0} |T\varphi_{j^* n_0} \gamma\|^\alpha}
$$
If $T$ is chosen so that
\begin{equation}
\label{FDer}
\| d\varphi_{j^* n_0, (j^*+1) n_0}|T\varphi_{j^* n_0} \gamma\|\geq
\left(R  e^{\epsilon(T-j^* n_0)} \right)^{-1},
\end{equation}
then the last inequality becomes
\begin{equation}
\label{APrioriCu}
K_{j^*+1}\leq
\frac{8 K_{j^*} \|d\varphi_{j^*n_0, (j^*+1) n_0}\|}
{\|d\varphi_{j^*n_0, (j^*+1) n_0}|T\varphi_{j^* n_0}\gamma\|^{1+\alpha}}
+ 2 R^2 e^{-(\lb_1' n_0(1-\alpha)-2 \epsilon)(T-j^* n_0)}
\end{equation}
Let us summarize what we have learned so far.

\begin{lemma} \label{HaveCu} For $n_0$ as above suppose that $T$ is
such that for every $j$ such that $j n_0\leq T$ estimates
(\ref{ScDer})--(\ref{FDer}) hold true and also the solution of
\begin{equation}
\label{DefK*}
\brK_{j+1}=
\frac{4 \brK_j \|d\varphi_{jn_0, (j+1) n_0}\|}
{\|d\varphi_{jn_0, (j+1) n_0}|T\varphi_{j n_0}\gamma\|^{1+\alpha}}
+ 2 R^2,\quad \brK_0=\brK
\end{equation}
satisfies
\begin{equation}
\label{SlExpGr}
\brK_j\leq \bar K e^{(T-j n_0)\epsilon}
\end{equation}
then inequality (\ref{SlGrCu}) holds.
\end{lemma}

Now we want to show that the set of points where either
(\ref{ScDer})--(\ref{FDer}) or (\ref{SlExpGr}) fail has density
less than $\varepsilon T$ except on a set of exponentially small
probability. The result for (\ref{SlExpGr}) follows from
Proposition \ref{Drift} applied to $\ln\brK_j.$ To see that
the conditions of this proposition are satisfied if $\alpha$ is
sufficiently small
it is enough to verify
that $\ln\brK_j$ has uniform drift to the left.

By Carverhill's extension of Oseledets' Theorem \cite{Cv} for every
point $x$ on $M$ and every unit vector $v$ in $T_xM$
\begin{equation}
\label{GrVec}
\frac{1}{n_0}\ \EXP\ \ln \| d\varphi_{n_0}(x) v\|\to\lb_1
\end{equation}
uniformly as $n_0\to\infty$ and \cite{BS} provides
exponential estimate for probabilities of large deviations.
Since
$$\|d\varphi_{n_0} (x)\|\leq \sum_{j=1}^N
\|d\varphi_{n_0} (x)v_j \| $$
where $\{v_j\}_{j=1}^N$ is any orthonormal frame,
the above mentioned results of \cite{BS} imply
that
$$
\frac{1}{n_0}\ \EXP\ \ln|\| d\varphi_{n_0}(x)\| \to\lb_1
\ \textup{as} \ \ \ n_0 \to \infty
$$
with exponential bound for large deviations. Thus
Proposition \ref{Drift} from the Appendix applies to $\ln \brK_j$
for a large enough $n_0$.

The fact that (\ref{ScDer})--(\ref{FDer}) fail rarely if $n_0$ is
sufficiently large and $r$ is sufficiently small follows from
Lemma \ref{SmDen}.

\subsection{Hyperbolic moments. Control of expansion.}
\label{hypmoments}

We now define the stopping time $\tau$ as the first moment when
(\ref{SlGrCu}), (\ref{ScDer}), and (\ref{FDer}) are satisfied
as well as
\begin{equation}
\label{Expansion}
\|d\varphi_\tau |\ T\gamma\|(x)\geq 1000
\end{equation}
and for each positive integer $j\leq\tau/n_0$ and some constant
$0<\tilde \lb_1<\lb_1$ we have
\begin{equation}
\label{BackDec}
\|d\varphi_{\tau, \tau-j n_0} |\ T\varphi_\tau\gamma\|\leq
e^{-\tilde \lb_1 j n_0}.
\end{equation}
Then the large deviations estimates of \cite{BS} guarantee that
property (\ref{Expansion}) has density close to 1.

\begin{lemma} For any $\varepsilon>0$ and any $0<\tilde \lb_1<\lb_1$
there exists a positive integer $n_0$ such that with probability
exponentially approaching to $1$ the fraction of integers $\tau$,
divisible by $n_0$, with the linearization
$d\varphi_{\tau, \tau-j n_0}|T \varphi_\tau \gamma$ contracting
exponentially backward in time for all integer $j$ between $0$ and
$\tau/n_0$ tends to $1$. More precisely,
$$
\PROB\left\{\frac{\#\{S\leq L: \forall 0\leq j \leq S,\ \tau=Sn_0\ \
\|d\varphi_{\tau, \tau-j n_0}|T \varphi_\tau \gamma\|
\leq e^{-\tilde \lb_1 j n_0}\}}{L}
\leq 1-\varepsilon \right\}
$$
decays exponentially in $L.$
\end{lemma}
\proof We first show how to prove a weaker statement with
``$\exists\varepsilon$'' instead of ``$\forall\varepsilon$''
(which is enough to prove Theorem \ref{ThHypT}) and then explain
briefly the changes needed to prove the sharp result.

Let $\tau_1$ be the first moment such that for each integer
$j\leq \frac{\tau}{n_0}$
\begin{equation}
\label{FHT}
\|d\varphi_{\tau_1, \tau_1-jn_0}|T\varphi_\tau\gamma\|
\leq e^{-\tilde \lb_1 j n_0}.
\end{equation}
We claim that $\tau_1$ has exponential tail. Indeed, let
\beq \nonumber
\beal
Y_j=Y_j(\theta) & =\left(\|d\varphi_{j n_0}|T\gamma\|(x)\
e^{-(\tilde \lb_1+\varepsilon)j n_0}\right)^\theta,\quad Y_0=1, \\
\textup{and} \quad
Z_j & =\|d\varphi_{j n_0}|T\gamma\|(x)\ e^{-\tilde \lb_1 j n_0},
\quad Z_0=1.
\enal
\eneq
Then \cite{BS} shows that if $n_0$ is sufficiently large and
$\varepsilon,  \theta$ are sufficiently small, then $Y_j$ is
a submartingale. Thus the first moment $\hj$ such that $Z_\hj>10$
has exponential tail. But there is at least one maximum $\brj$ of
$Z_j$ between $0$ and $\hj.$ Then $\brj$ satisfies (\ref{FHT}).

Now define $\tau_k$ inductively so that $\tau_{k+1}>\tau_k$ is the first
moment such that for every $j\leq \frac{\tau_{k+1}-\tau_k}{n_0}$
$$
\|d\varphi_{\tau_{k+1}, \tau_{k+1}-jn_0}|T\varphi_{\tau_{k+1}}\gamma\|
\leq e^{-\tilde \lb_1 n_0 j} .
$$
Then $\tau_{k+1}-\tau_k$ have exponential tails, so by Lemma \ref{MART}
there exists $c$ such that $ \PROB\{\frac{\tau_k}{k}\geq C\}$ decays
exponentially in $k.$ However all $\tau_k$ satisfy (\ref{FHT}).
This proves the result with $\varepsilon=1-\frac{1}{C}.$ To get
the optimal result one should note that $\PROB\{\tau_1=n_0\}\to 1$ as
$n_0\to\infty$ and apply the arguments of Lemma \ref{MART}. We leave
the details to the reader. \qed

Now we want to verify conditions (b), (c), and (d) of Theorem \ref{ThHypT}
with $\brgamma_r$ replaced
by $\hgamma=\varphi_\tau B_{\tau,\tau}(x).$ Once we prove this we get
from (c) that the main restriction on $B_{\tau,\tau}(x)$ is for
$k=\tau$ so that $\brgamma_r=\hgamma$ and then (a) will also be true.
Now (b) is true by Lemma \ref{HaveCu}. We will establish (c) and (d)
by induction. Namely we suppose that (c) is true for $k\geq k_0.$
Then for every $y \in \hgamma$
\beq \nonumber
\beal
\left| \ln\|d\varphi_{\tau, \tau-(k_0+1) n_0}|T\hgamma \|(x) \right. &
\left. -
\ln\|d\varphi_{\tau, \tau-(k_0+1) n_0}|T\hgamma \|(y) \right|\leq \\
\sum_{m=0}^{k_0} \left|
\ln\|d\varphi_{\tau-mn_0, \tau-(m+1) n_0}|T\hgamma\|(x) \right. &
\left. -
\ln\|d\varphi_{\tau-mn_0, \tau-(m+1) n_0}|T\hgamma\|(y)\right|\leq \\
\sum_{m=0}^{k_0} \left|
\ln\|d\varphi_{\tau-mn_0, \tau-(m+1) n_0}|T\hgamma\|(x) \right. &
\left. -
\ln\|d\varphi_{\tau-mn_0, \tau-(m+1) n_0}|T\hgamma\|(y)\right|+ \\
\sum_{m=0}^{k_0} \left|
\ln\|d\varphi_{\tau-mn_0, \tau-(m+1) n_0}|T\hgamma\|(y) \right. &
\left. -
\ln\|d\varphi_{\tau-mn_0, \tau-(m+1) n_0}|T\hgamma\|(y)\right|.
\enal
\eneq
Denote the left term by $I$ and the right term by $\RmII$
respectively. Now by (\ref{SlGrCu}) and (\ref{ScDer})
\beq
\beal
I\leq\sum_{m=0}^{k_0} R e^{\epsilon m} K e^{\epsilon m}
\rho^\alpha(\varphi_{\tau,\tau-m n_0} x, \varphi_{\tau,\tau-m n_0} y)
\leq \\
\sum_{m=0}^{k_0} K R e^{-(\lb_1' \alpha n_0-2\epsilon) m}
\rho^\alpha(x,y)\leq \Const r^\alpha.
\enal
\eneq
On the other hand
\beq
\beal
\RmII \leq
\sum_{m=0}^{k_0} \frac{\|d^2\varphi_{\tau-mn_0,\tau-(m+1)n_0}\|}
{\|d\varphi_{\tau-mn_0, \tau-(m+1)n_0}\|} & \
\rho(\varphi_{\tau,\tau-m n_0} x, \varphi_{\tau,\tau-m n_0} y)\leq
\\
\sum_{m=0}^{k_0} R^2 e^{-(\lb_1' n_0-2\epsilon) m} &
\ \rho(x,y)\leq \Const r.
\enal
\eneq
Hence (\ref{ScDer}) and (\ref{FDer})
\begin{equation}
\label{Distortion}
\left| \ln\|d\varphi_{\tau, \tau-(k_0+1) n_0}|T\hgamma\|(x)-
\ln\|d\varphi_{\tau, \tau-(k_0+1) n_0}|T\hgamma\|(y) \right|\leq
C(R) \rho(y_1, y_2).
\end{equation}
Thus for all $y$
\begin{equation}
\label{GrDer}
\|d\varphi_{\tau-(k_0+1) n_0,\tau} |
T\varphi_{\tau-(k_0+1) n_0}\gamma\|(y)
\geq
\exp\left(\tilde \lb_1 k n_0-C(R) r\right)\geq
\exp\left(\lb_1' k n_0\right)
\end{equation}
if $\tilde \lb_1-\lb_1'\geq C(R) r.$
(\ref{GrDer}) implies that (c) is valid for $k_0-1.$ Thus, we obtain
(c) for all $k.$ Now repeating the proof of (\ref{Distortion}) with
$x$ and $y$ replaced by $y_1$ and $y_2$ (and using (\ref{GrDer})
instead of (\ref{BackDec})) we obtain (d).
This completes the proof of Theorem \ref{ThHypT}. \qed

\begin{remark}
The term {\it hyperbolic time} was introduced in \cite{A} but the notion
itself was used before, e.g. in \cite{P1, P2, J, Y}. Considerations of
this section are similar to \cite{ABV, D2} but the additional difficulty
is that in those papers the analogue of (\ref{SlGrCu}) was true by
the general theory of partially hyperbolic systems \cite{HPS} whereas
here additional arguments in spirit of \cite{P1, P2} were needed
to establish it.

One interesting question is how large can $\alpha$ be so that
Theorem \ref{ThHypT} still holds.
We note that $\alpha$ appears in (\ref{APrioriCu}) twice. So we want $\alpha$
to be as large as possible to control the first part and we want $\alpha$
to be small to control the second term. In general, the optimal choice of
$\alpha$ should depend on the ratio of leading exponents.
We refer to \cite{CL, L, PSW, JPL}
for the discussion of this question.
\end{remark}

\section{Construction of the partition.}
\label{ScPart}
We are now ready to describe a partition
$\gm=\bigcup_{j\in \Z_+}\gm(j).$ It will be defined inductively.
Each of $\gm(j)$'s is a finite union of intervals. As $j$ tends to
infinity size of intervals tends to zero and they fill up $\gm$.
To simplify the notation we assume that Theorem \ref{ThHypT} is true
with $n_0=1$. This can be achieved by rescaling the time. Fix
an orientation from left to right on $\gamma$.

Suppose $\gamma(1), \gamma(2), \dots, \gamma(m)$ are already defined
in an $\cF_m$-measurable way. Let
\beq
K_{m+1}=\{x\in \gm:\ \tau(x)=m+1\}.
\eneq
By definition $K_{m+1}$ is a finite union of intervals.
Let $U_{m+1}=\varphi_{m+1} K_{m+1}.$ We call {\it an obstacle}
any point on the boundary of either $K_{m+1},$
$\bigcup_{j=1}^{m} \gamma(j)$ or $\gamma.$ Fix $r$ satisfying
Theorem \ref{ThHypT}. Let $C$ be a connected component of $U_{m+1}$
and $a$ and $b$ be its left and right endpoints with respect to left-right
orientation induced by $\varphi_{m+1}$. If distance from $b$ to
the closest image of an obstacle to the right on $\varphi_{m+1}(\gm)$
is less than $\frac{r}{2}$ and $b'$ is this image, then put $\tb=b'$.
Otherwise let $\tb$ be a point at distance $\frac{r}{100}$ from $b.$
Define $\ta$ similarly. Consider the set $W_{m+1}=\bigcup_C \ta\tb.$
Divide $W_{m+1}$ into the segments of
lengths between $\frac{r}{100}$ and $\frac{r}{50}$ and denote this
partition by $V_{m+1}$. Now we define partition of a subset of
$\gm \setminus \bigcup_{j=1}^{m} \gamma(j)$ by pulling back  along
$\varphi^{-1}_{m+1}$ the partition $V_{m+1}$
\beq
\gm(m+1)= \varphi^{-1}_{m+1} V_{m+1}.
\eneq
To justify that this algorithm produces a partition which covers all
of $K_{m+1}$ we need to check that length of each component is at least
$\frac{r}{100}.$ To do this we argue by contradiction. Otherwise,
there would be two obstacles $x', x''$ neither of which is from $K_{m+1}$
such that $ \rho(\varphi_{m+1}x', \varphi_{m+1} x'')\leq \frac{r}{100}$
and a point from $U_{m+1}$ between them. At least one of the obstacles
would have to come from $\bigcup_{j=1}^m \gamma(j).$ Let $x'$ be such
an obstacle. Since both points are close to $U_m$ for each $n\leq m+1$
we have
$$
\rho(\varphi_n x', \varphi_n x'')\leq\frac{r}{100}\ e^{-\lb_1'(m-n)}.
$$
But in this case the interval $[x',x'']$ in $\gm$ with endpoints $x'$ and
$x''$ would be added to our partition at a previous step of the algorithm.

Denote by $\gamma=\bigcup_{j\in \Z_+}\gamma_j$ the partition which is
made out of the partition $\gamma=\bigcup_{j\in \Z_+}\gamma(j)$ by
renumerating intervals of this partition in length decreasing order.
Let us summarize the outcome.
\begin{proposition}
\label{PrPart}
We can partition $\gamma=\bigcup_{j\in \Z_+}\gamma_j$ in such a way
that

\noindent
(a) there exists a positive integer $n_j$ such that
$\|d\varphi_{n_j}|T\gamma\| \geq 100$ and length
$l(\varphi_{n_j} \gamma_j)\geq  \frac{r}{100}$;

\noindent
(b) for each positive integer $m\leq n_j$ and lengths of
the corresponding curves we have  $l(\varphi_m \gamma_j)\leq
l(\varphi_{n_j} \gamma_j) e^{-\lb_1'(n_j-m)};$

\noindent (c) \ $\left| \ln \|d\varphi_{n_j}| T\gamma\|(x')-
\ln \|d\varphi_{n_j}| T\gamma\|(x'')\right|\leq
\Const \rho^\alpha(\varphi_{n_j} x', \varphi_{n_j} x'')$
for every pair $x', x''\in \gamma_j$;

\noindent (d)\  for some $\alpha>0$ and each pair
$x', x''\in \gamma_j$ we have
$\left|v(x', n_j)-v(x'',n_j)\right|\leq \Const
\rho^{\alpha} (\varphi_{n_j} x', \varphi_{n_j} x'')$ ,
where $v(x,n)$ denote the unit tangent vector to
$\varphi_n\gamma$ at $\varphi_n (x);$

\noindent
(e) Let $j(x)$ be such that $x\in\gamma_{j(x)}.$ Then
$\EXP\ n_{j(x)}\leq \Const$ and\ $\PROB\{n_{j(x)}>T\}
\leq C_1 e^{-C_2 T}$ for some positive $C_1, C_2$ and any
$T>0$;
\end{proposition}

This Proposition is designed to allow application of
Theorem \ref{ThHypT} so that we can use regularity and geometric
properties of $\gm_j$'s at stopping times $\tau_j$'s.

\section {Construction of the measure with almost sure nonzero drift.}
\label{ScDefMes}
Now we construct a random Cantor set $I\subset \gamma$ and
a probability measure $\mu$ supported on $I$ such that $\mu$--almost all
points have a nonzero drift. This construction goes along the same line
with the construction in Section \ref{drift} of the Cantor set $I_\infty$
in the unit interval and a probability measure $\mu_\infty$ on $I$
such that $\mu_\infty$--almost all points have nonzero drift.

Choose a direction $\ve\in\R^N.$ Let $\theta$ be a small parameter
which we let to zero in the next section. We say that a curve is
{\it $\ve$-monotone} if its projection to $\ve$ is monotone. Now
we describe construction of a Cantor set $I\subset \gamma$ and
a probability measure $\mu$ on $I$ by induction.
This Cantor set $I$ at $k$-th step of induction consists of
countable number of segments numerated by $k$-tuples of positive
integers.

Denote $k$-tuples $(j_1,\cdots,j_k)\in \Z_+^k$ and
$(n_1,\cdots,n_k)\in \Z_+^k$ by $J_k$ and $N_k$ respectively.
Let $|N_k|=\sum_{j=1}^k n_j$.

The first step of induction goes as follows. Let $\gamma_j,$ $n_j$
be the sequence of pairs: a curve and an integer, described in
Proposition \ref{PrPart}. Let $\theta$ be a small positive number.
If $\varphi_{n_j} \gm_j$ is $\ve$-monotone put $\sigma(j)$ equal
$\varphi_{n_j} \gm_j$ without the segment of length $\theta r$, which
we cut off from the $\ve$-bottom point of $\varphi_{n_j}\gm_j.$
Otherwise $\sigma(j)=\varphi_{n_{j}} \gamma_j$ with no cut off.
Let $\gamma(J_1)=\varphi_{n_j}^{-1}\sigma(j)$ and
$N_1(J_1)=n_{j}$ for $J_1=j$

Suppose a collection of disjoint segments
$\{\gamma(J_k)\}_{J_k \in \Z_+^k} \subset \gm$ is defined as above
and multiindices $N_k$ (resp. $J_k)$ are defined as the corresponding
set of hyperbolic times multiindexed by $J_k$ segments. Then
\beq \label{partition}
I_k=\cup_{J_k \in \Z_+^k} \gamma(J_k) \subset I_{k-1}\subset \dots
\subset I_1 \subset \gm
\eneq
is the $k$-th order of construction of the random Cantor set $I$
(cf. with an open set $I_k$ from Section \ref{drift}).

The $(k+1)$-st step goes as follows. Pick a segment $\gm(J_k)$ of
partition (\ref{partition}). Consider the partition of the curve
\beq
\varphi_{|N_{J_k}|} \gamma(J_k)=
\bigcup_{j_{k+1}\in \Z_+} \tilde \gamma(J_k,j_{k+1})
\eneq
defined in Section \ref{ScPart} and let $n_{J_k,j_{k+1}}$ be
the corresponding hyperbolic times for $\tilde \gamma(J_k,j_{k+1})$
from Proposition \ref{PrPart}. For brevity denote $|N_k(J_k)|$
by $n^{(k)}$ and $|N_k(J_k)|+n_{(J_k,j_{k+1})}$ by $n^{(k+1)}$.
If the curve $\varphi_{n^{(k)},n^{(k+1)}} \tilde\gamma_{(J_k,j_{k+1})}$
is $\ve$-monotone we let $\sigma(J_k,j_{k+1})$ be
$\varphi_{n^{(k)},n^{(k+1)}} \tilde\gamma_{(J_k,j_{k+1})}$ with
cut off of the segment of length $\th$ starting from the $\ve$-bottom.
Otherwise, $\sigma(J_k,j_{k+1})$ equal
$\varphi_{n^{(k)},n^{(k+1)}}\tilde\gamma(J_k,j_{k+1})$ with no cut off.
Then a segment
\beq
\gamma(J_k,j_{k+1}) = \varphi^{-1}_{n^{(k+1)}}
\sigma(J_k,j_{k+1})
\eneq
with $j_{k+1}\in Z_+$ this defines the $(k+1)$-st order partition
$\{\gamma(J_{k+1})\}_{J_{k+1} \in \Z_+^{k+1}} \subset \gm$ and
the $k$-order set $I_{k+1}=
\cup_{J_{k+1} \in \Z_+^{k+1}} \gamma(J_{k+1}) \subset \gm$.

We now describe a sequence of measures $\mu_k$'s on
$I_k \subset \gamma$ with $k\in\Z_+$ respectively.
Let $\mu_0$ be the arclength on $\gamma.$ Suppose $\mu_k$ is already
defined on $I_k$.
Consider $\{\gamma(J_{k+1})\}_{J_{k+1}\in \Z_+^{k+1}}.$
If $\varphi_{n^{(k+1)}}\gamma(J_{k+1})$ is {\it not}
$\ve$-monotone we let
$\mu_{k+1}|_{\gamma(J_{k+1})}=\mu_{k}|_{\gamma(J_{k+1})}.$
Otherwise, $\mu_{k+1}|_{\gamma(J_{k+1})}=\rho_{jk}
\mu_{k}|_{\gamma(J_{k+1})}$, where $\rho_{jk}$ is a normalizing
constant.

\begin{lemma}
\label{Mono} Let $k$ be an integer. If $r$ is sufficiently small and
$\gamma \subset \R^N$ is $(K,\alpha)$-smooth as in Theorem \ref{ThHypT},
then if we consider partition of $\gamma$ up to order $k+1$, then
for each multiindex $J_k\in \Z^k_+$ the corresponding $k$-th order
curve $\gamma(J_k)\subset \gamma$ satisfy the property: for
any positive integer $j_{k+1}$ the $(k+1)$-st order curve
$\gamma(J_{k+1}) \subset \gamma(J_k)$ has $\ve$-monotone with
positive probability, i.e.
$$
\PROB\{\varphi_{n^{(k+1)}}\gamma(J_{k+1})\ \
\text{is}\ \ \ve-\text{monotone}\ |\ \ \cF_{n^{(k)},n^{(k+1)}}\} >c
$$
for some positive $c$ and $c$ is uniform for
all $(K,\alpha)$-smooth curves.
\end{lemma}

\proof Pick a point $x \in \gamma(J_{k+1})$. By assumption (D) of
hypoellipticity on the unit tangent bundle $SM$ for the flow (\ref{SF})
probability that the angle between $\ve$ and
$T\varphi_{n^{k+1}}\gamma(x)$ makes  less than $1^\circ$ is positive.
By definition $\varphi_{n^{(k+1)}}\gamma(J_{k+1})$ is
$(K,\alpha)$-smooth. Thus if $r$ is small enough, then the tangent
vectors to $\varphi_{n^{k+1}}\gamma$ are close to
$T\varphi_{n^{k+1}}\gamma (x)$ with large probability, where
$x$ is a point on $\gamma(J_{k+1})$. This completes the proof. \qed

Recall that $\theta>0$ is a fraction of $\varphi_{n_j}\gamma_j$
we cut off from $\varphi_{n_j}\gamma_j$ on the $j$-th step,
provided $\varphi_{n_j}\gamma_j$ is $\ve$-monotone.
Let $\mu=\mu(\theta)$ denote the weak limit of $\mu_k$'s
$$
\mu=\lim_{k\to \infty} \mu_k.
$$
\begin{lemma}
For almost every realization of the Brownian motion $\{\th(t)\}_{t\geq 0}$
and $\mu$--almost every $x$
$$
\liminf_{t\to\infty} \frac{\langle x_t, e\rangle}{t}>0.
$$
\end{lemma}

\proof The first step is to show that for any $s$ for almost all
realizations of the Brownian motion $\{\th(t)\}_{t\geq 0}$
\begin{equation}
\label{DiscrTmDrift}
\liminf_{t\to\infty}
\frac{\langle x_{n^{(k)}}^{(k)}, e\rangle}{n^{(k)}}>0
\end{equation}
Applying Proposition \ref{PrPart} (e) and Lemma \ref{MART} we get
that there exists a constant $C>0$ such that
$$
\limsup_{k\to\infty} \frac{n^{(k)}}{k}<C
$$
almost surely. Therefore, to prove (\ref{DiscrTmDrift})
it suffices to show that
\begin{equation}
\label{TmChDrift}
\liminf_{k\to\infty}\frac{\langle x_{n^{(k)}}^{(k)}, e\rangle}{k}>0
\end{equation}
However by Lemma \ref{Mono} there exists $c$ such that
$\EXP \langle x_{n^{(k+1)}}^{(k+1)}-x_{n^{(k)}}^{(k)},e\rangle >c$
uniformly in $k,s.$  (This is because
$\EXP \langle x_{n^{(k+1)}}^{(k)}-x_{n^{(k)}}^{(k)},e\rangle =0$
and
$$
\langle x_{n^{(k+1)}}^{(k)},e \rangle -
\langle x_{n^{(k)}}^{(k)},e \rangle\ \geq 0
$$
with
strict inequality having positive probability by Lemma \ref{Mono}.)
Hence (\ref{TmChDrift}) follows by Lemma \ref{MART}. Therefore
(\ref{DiscrTmDrift}) is established. \qed

Now we apply the following estimate.
\begin{lemma} (\cite{CSS2}, Theorem 1)
Let
$$
\Phi_{s,t}=\sup_{s\leq\tau\leq t} |x_\tau-x_s|, \quad
\tPhi_{s,t}=\frac{\Phi_{s,t}}{\max(1,t-s)} $$
then there exists a constant $C$ such that for all $s$ and $t$
$$
\EXP\left(\exp\left\{\frac{\tPhi_{s,t}^2}{\max(1,\ln^3 \tPhi_{s,t})}
\right\}\right)<C.
$$
\end{lemma}

Combining this lemma with Proposition \ref{PrPart} (e) we obtain
that there are positive constants $\alpha$ and $D$ such that
$$
\EXP\left(\exp\left\{\alpha\sup_{n^{(k)}<\tau<n^{(k+1)}}
\left|x_\tau(s)-x_{n^{(k)}}\right|\right\}\right)<D.
$$
Using Borel-Cantelli's lemma we derive from this that almost surely
$$
\lim\sup_{k\to\infty}\frac{\sup_{n^{(k)}<\tau<n^{(k+1)}}
\left|x_\tau(s)-x_{n^{(k)}}\right|}{\ln k}<+\infty .
$$
Therefore for any $s$ and for almost all realizations
of the Brownian motion $\{\th(t)\}_{t\geq 0}$ we have
$$
\liminf_{\tau\to\infty} \frac{\langle x_\tau(s),e\rangle}{\tau}>0 .
$$
By Fubini Theorem we have that for almost every realization of
the Brownian motion $\{\th(t)\}_{t\geq 0}$ the set
$$
\left\{ s: \quad \liminf_{\tau\to\infty}
\frac{\langle x_\tau(s),e\rangle}{\tau}>0 \right\}
$$
has full measure. \qed

\section{Hausdorff dimension of $\mu.$} \label{ScHD}

In this section we complete the proof of Theorem \ref{HDC}
by establishing the following fact. Recall that $\theta>0$ is
a fraction of $\varphi_{n_j}\gamma_j$ we cut off from
$\varphi_{n_j}\gamma_j$ on the $j$-th step, provided
$\varphi_{n_j}\gamma_j$ is $\ve$-monotone. Consider
the measure $\mu$ we constructed in the previous Section.

\begin{proposition}
\label{PrSmallMes} With notations above we have that
as $\theta\to 0$ Hausdorff dimension of the measure
$\mu=\mu(\theta)$ tends to $1$:\ $\HD(\mu(\theta))\to 1.$
\end{proposition}

Let us recall the following standard principle.

\begin{lemma} [Mass distribution principle]
Let $S$ be a compact subset of a Euclidean (or metric) space such
that there exists a probability measure $\nu$ such that $\nu(S)=1$
and for each $x$ we have $\nu(\cB(x,r))\leq C r^s$
for some positive $C$ and $s$. Then $\HD(S)\geq s.$
\end{lemma}
Proposition \ref{PrSmallMes} is a direct consequence of
the following statements.

\begin{lemma}
\label{MesDist}
Let $\gamma$ be a smooth curve in $\R^N$. Suppose there exist
a nested sequence of partitions
\beq
\gamma \supset \bigcup_{J_1\in \Z^1_+} \gamma(J_1) \supset \dots
\supset \bigcup_{J_k\in \Z^k_+} \gamma(J_k) \supset \dots
\eneq
and probability measures $\mu_0, \mu_1, \dots ,\mu_k, \dots$
supported on $\gamma, \bigcup_{J_1\in \Z^1_+} \gamma(J_1), \dots,$
$\bigcup_{J_k\in \Z^k_+} \gamma(J_k), \dots $ respectively such
that $\mu_0$ is the normalized arclength on $\gm$ and so on $\mu_k$
is the normalized arclength on $\bigcup_{J_k\in \Z^k_+} \gamma(J_k)$.
Then if we have

\noindent (b) for all $J_k \in \Z^k_+$ length of the corresponding
interval $\gamma(J_k)$ is bounded by $l(\gamma(J_k))\leq 100^{-k};$

\noindent (b) for each $l>k$ we have
$\mu_l(\gamma(J_k))=\mu_k(\gamma(J_k));$

\noindent (c) $\frac{d\mu_{k+1}}{d\mu_k}(x)\leq \rho_k(x)$ for every
point $x \in \bigcup_{J_k\in \Z^k_+} \gamma(J_k)$,
where $\rho_k<1+\delta.$

Let $\mu=\lim_{k\to\infty} \mu_k$ in the sense of weak limit. Then
$\HD(\mu)\geq d(\delta)$, where $d(\delta)\to 1$ as $\delta\to 0.$
\end{lemma}

\begin{lemma}
\label{LongTimeDistortion}
For each $\delta>0$ there exists $\theta>0$ such that the densities
of $\frac{d\mu_{k+1}}{d\mu_k}(x)$ used to define measures
$\mu_{k+1}$ knowing $\mu_k$ satisfy condition  (c)
of Lemma \ref{MesDist}.
\end{lemma}
{\it Proof of Lemma \ref{MesDist}.} We prove that for any segment $I$
we have $\mu(I)\leq\Const|I|^{1-\beta},$ where $\beta\to 0$
as $\delta\to 0.$ Let $k(I)=\frac{|\ln|I\|}{\ln 100}$ and
$a$ and $b$ be the endpoints of $I.$ Let  $\ta$ be the left
endpoint of the $k$-th partition containing $a$ and
$\tb$ be the right endpoint of the $k$-th partition
containing $b.$
Then
\begin{equation}
\label{SmallMeasure}
\mu(I)\leq \mu([\ta,\tb])=\mu_k([\ta, \tb])\leq
(1+\delta)^k \mu_0([\ta, \tb])\leq 3(1+\delta)^k |I|
\leq 3|I|^{1-\beta}
\end{equation}
where $\beta=\frac{\ln(1+\delta)}{\ln 100}.$ Thus,
$\beta\to 0$ as $\dt\to 0.$ Application of the mass distribution
principle implies that $\HD(\mu)\geq 1-\beta.$ \qed

{\it Proof of Lemma \ref{LongTimeDistortion}.}
Recall the notation of Section \ref{ScDefMes}.
We need to show that
\begin{equation}
\label{SmallCut}
\sup_{J_{k+1}}
\frac{|\varphi^{-1}_{n^{(k+1)}} \sigma(J_{k+1})|}
{|\varphi^{-1}_{n^{(k)}} \gamma(J_{k+1})|}
\to 1, \quad \theta\to 0
\end{equation}
By construction
$$
\frac{|\varphi_{n^{(k+1)}, n^{(k)}} \sigma(J_{k+1})|}
{|\gamma(J_{k+1})|}\geq 1-\left(\frac{\theta}{r/100}\right) .
$$
Hence to prove (\ref{SmallCut}) it is enough to show that there is
a constant $C$ independent of $j, k, l$ such that for any interval
$I\subset\gamma(J_{k+1})$
$$
\frac{|\varphi^{-1}_{n_j^{(k)}} I|}{|\varphi^{-1}_{n_j^{(k)}}
\gamma(J_{k+1})|}\leq C\frac{|I|}{|\gamma(J_{k+1})|} .
$$
To do so it is enough to show that there is a constant $\brC$
such that  for every pair $y_1, y_2\in \gamma(J_{k+1})$
$$
\frac{\| d\varphi_{n_j^{(k)}}^{-1}| T\gamma(J_{k+1})\|(y_1)}
{\| d\varphi_{n_j^{(k)}}^{-1}| T\gamma(J_{k+1})\|(y_2)}\leq\brC.
$$
But by Proposition \ref{PrPart} there are constants $C_1, C_2,$ and
$C_3$ such that
\beq \nonumber
\beal
\left|\ln\|d\varphi_{n_j^{(k)}}^{-1}| T\gamma(J_{k+1})\|(y_1)-
\right. & \left.
\ln\|d\varphi_{n_j^{(k)}}^{-1}| T\gamma(J_{k+1})\|(y_2)\right|\leq \\
\sum_{m=1}^k
\left|\ln\|d\varphi_{n^{(m)}, n^{(m-1)}}| \right. & \left.
T\varphi_{n^{(m)}}\|(\varphi_{n^{(k)}, n^{(m)}} y_1) \right. \\
- \left.
\ln\|d\varphi_{n^{(m)}, n^{(m-1)}}|
T \right. & \left.\varphi_{n^{(m)}}\|
(\varphi_{n^{(k)}, n^{(m)}} (y_2))\right| \leq  \\
C_1\ \sum_m \rho^\alpha(\varphi_{n^{(k)}, n^{(m)}} y_1,
\varphi_{n^{(k)}, n^{(m)}} (y_2))\leq
C_2\ & \sum_m 100^{(m-k)\alpha} \rho^\alpha(y_1, y_2)
\leq C_3\ r^\alpha.
\enal
\eneq
This completes the proof. \qed

\section{Proof of Theorem \ref{FullDim}.}
\label{fullHD}

Let $\cG$ denote the foliation of $\T^N$ by curves
$$
\{x^1=c^1, x^2=c^2\dots x^{N-1}=c^{N-1}\}.
$$
By (\ref{SmallMeasure}) for each $\beta>0$ and each leaf $\gamma_c$
of $\cG$ almost surely there exists a measure $\mu_c$ on $\gamma_c$
such that $\mu_c(I)\leq 3 |I|^{1-\beta}$ and $\mu_c(\bL_\theta)=1.$
Let $\mu=\int \mu_c dc.$ Then by Fubini Theorem almost surely for
any cube $\cC$ of side $r$ we have $\mu(\cC)\leq 3 r^{N-\beta}$ and
$\mu(\bL_\theta)=1.$ The application of the mass distribution
principle completes the proof. \qed

\begin{appendix}
\section{Large deviations.}
Here we collect some estimates used throughout the proof of
Theorem \ref{HDC}.
\begin{lemma}
\label{MART} Let $\cF_j$ be a filtration of $\sigma$-algebras and
$\{\xi_j, \}$ be a sequence of $\cF_j$-measurable
random variables such that

\noindent (a) there exist $C_1, \lb$ such that for every $|s|\leq
\lb$ we have $\EXP (  e^{s\xi_{j+1}} | \cF_j) \leq
C_1;$

\noindent (b) there exists $C_2$ such that  $\EXP (\xi_{j+1}
|\cF_j ) \leq C_2.$

Then  for each $\epsilon>0$ the probability $$ \PROB
\left\{\sum_{j=0}^{N-1} \xi_j \geq (C_2+\epsilon) N\right\} $$
decays exponentially in $N$.
\end{lemma}

\proof Consider
$$ \Phi_n(s)=\exp\left\{\left( \sum_{j=0}^{n-1}
\xi_j- \left(C_2+\frac{\epsilon}{2} \right)n\right)s\right\}. $$
Then (a) and
(b) imply that  $\Phi_n(s)$ is a supermartingale if $s$ is
sufficiently small. Hence $\EXP\Phi_n(s) \leq\EXP\Phi_0(s)=1,$ and
so $$ \EXP\ \exp\left\{\left( \sum_{j=0}^{n-1} \xi_j- (C_2+
\epsilon)n\right)s\right\}\leq \exp\left(-\frac{  n \epsilon s
}{2}\right)~, $$ which proves the lemma. \qed

\begin{lemma}
\label{SmDen}
Let $\cF_j$ be a filtration of $\sigma$-algebras and
$\{\xi_j, \}$ be a sequence of $\cF_j$-measurable
random variables such that there exists constant $C_1$ such that
\begin{equation}
\label{UniExpMom} \EXP(\xi_{j+1}|\cF_j)\leq C_1
\end{equation}
then for every $\epsilon, \varepsilon>0$ there is $R>0$ such that $$
\PROB\left\{\frac{\Card\{n\leq N: \xi_j\leq R e^{\epsilon (n-j)}~~
{\rm for ~all} ~ 0 \leq j < n  \}}{N} \leq 1-\varepsilon\right\}
$$ tends to zero exponentially fast in $N$.
\end{lemma}

\proof
We say that a pair $(j,n)$ is $R$-bad if
$$\xi_j>R e^{\epsilon (n-j)} . $$
By (\ref{UniExpMom})
\begin{equation}
\label{BadUnlikely}
\PROB\{(j,n) \;\text{is $R$-bad}\}\leq
\frac{C_1}{R} e^{-\epsilon(n-j)}.
\end{equation}
Now given $k$ let $B_{R}(k)$ be the number of $n>k$ such that
$(k, n)$ is $R$-bad. By (\ref{BadUnlikely})
$$
\EXP(B_{R}(k+1)|\cF_k)\leq
\frac{C_2 e^{-\epsilon}}{R(1-e^{-\epsilon})} \to 0
$$
as $R\to\infty.$ Thus by Lemma \ref{MART} there exists $R$
such that
$$
\PROB\left\{\sum_{k=1}^N B_{R}(k)\geq \varepsilon N\right\}
$$
decays exponentially in $N$. This completes the proof
of the lemma. \qed

\begin{proposition}
\label{Drift}
Let $x_j$ be (in general, non-homogeneous) random walk on
$\integers.$ Suppose that there exist constants $C_1, C_2, C_3$
such that

\noindent (a) there exist $m$ such that for every $x_j>m$ we have
$\EXP(x_{j+1}-x_j|x_j)\leq -C_1;$

\noindent (b) for every $x_j$ and every $\zeta<C_2$ we have
$\EXP(e^{\zeta(x_{j+1}-x_j)}|x_j)\leq C_3. $

Fix $\delta>0.$ Let $F(M)$ denote the set of $j$ such that for all
$k<j$
$$
x_k\leq \max(x_j,m)+M+\delta(j-k).
$$
Then for every $\varepsilon>0$ there exists $M>0$ such that
$$
\PROB\left\{\frac{\Card\{F(M)\bigcap [1, N]\}}{N}\leq
1-\varepsilon \right\}
$$
decays exponentially in $N.$
\end{proposition}

\proof Let $\tau_1<\tau_2\dots<\tau_k<\dots$ be the consecutive
returns of $x_j$ to $\{x\leq m\}.$ Let
$$
t_j=\tau_{j+1}-\tau_j, \quad X_j=\max_{\tau_{j-1}<l<\tau_j} x_l.
$$
\begin{lemma}
$t_j$ and $X_j$ have exponential tails.
\end{lemma}

\proof It suffices to prove it for $t_1$ and $X_1$ and the assumption
that $x_0\leq m.$ Clearly it suffices to condition on $x_1>m$ since
otherwise $t_1=1,$ $X_1\leq m.$ Then (b) implies that for small
$\varepsilon_1, \varepsilon_2$
$$
y_j=e^{\varepsilon_1 x_j+\varepsilon_2 j} 1_{\{j\leq \tau_1\}}
$$
is a supermartingale. Thus
\begin{equation}
\label{NotFar}
\EXP\ y_j\leq \EXP\ y_1\leq C_4.
\end{equation}
On the other hand
$$
\EXP\ y_j\geq \PROB\{\tau_1>j\}\ e^{\varepsilon_1 m+\varepsilon_2 j}.
$$
Hence
$$
\PROB\{\tau_1>j\} \leq C_5 e^{-\varepsilon_2 j}
$$
where $C_5=C_4 e^{-\varepsilon_1 m}.$
Now
$$
\EXP\ e^{\varepsilon_1 X_1}\leq \EXP\left( \sum_{j=1}^{\tau_1}
e^{\varepsilon_1 x_j}\right)
\leq \sum_{j=1}^\infty \EXP e^{\varepsilon_1 x_j} \sum_{k=j}^{\infty}
\PROB\{\tau_1>k\}\leq
C_4 \sum_j \frac{C_5 e^{-\varepsilon_2 j}}{1-e^{-\varepsilon_2 j}}
<\infty .
$$
This completes the proof.  \qed

The rest of the proof of Proposition \ref{Drift} is similar to
the proof of Lemma \ref{SmDen}. We say that the pair $(k,j)$ is
bad if
$$
x_k>\max(x_j,m)+M+\delta(j-k).
$$
If $(k,j)$ is bad then
$$
j-k<\frac{x_k-m+M}{\delta}.
$$
Let $B_l(M)$ be the number of bad pairs $(k,j)$ such that
$\tau_{l-1}<k<\tau_l.$ By the previous lemma
$\EXP\ B_l(M)<\infty$ and so by dominated convergence theorem
$\EXP\ B_l(M)\to 0$ as $M\to \infty.$ Hence by Lemma \ref{MART}
the number of bad pairs such that $k<\tau_N$ is less than
$\varepsilon N$ except on a set of exponentially small probability.
Since $\tau_N\geq N$ the proposition follows. \qed

\end{appendix}

\end{document}